\documentclass[12pt,leqno]{amsart}  
\usepackage{amsmath,amstext,amsthm,amssymb}
\usepackage[colorlinks,pagebackref,hypertexnames=false]{hyperref} % citecolor=red,
\usepackage[backrefs]{amsrefs}

\allowdisplaybreaks

\newtheorem{lemma}[equation]{Lemma} 
\newtheorem{proposition}[equation]{Proposition} 
\newtheorem{theorem}[equation]{Theorem} 
\newtheorem{corollary}[equation]{Corollary}

\def\Label #1 {\label{#1}}

\def\norm#1.#2.{\lVert#1\rVert_{#2}}
\def\Norm#1.#2.{\bigl\lVert#1\bigr\rVert_{#2}}
\def\NOrm#1.#2.{\Bigl\lVert#1\Bigr\rVert_{#2}}
\def\NORm#1.#2.{\biggl\lVert#1\biggr\rVert_{#2}}
\def\NORM#1.#2.{\Biggl\lVert#1\Biggr\rVert_{#2}}

\def\ip#1,#2.{\langle #1,#2\rangle}
\def\Ip#1,#2.{\bigl\langle#1,#2\bigr\rangle}
\def\IP#1,#2.{\Bigl\langle#1,#2\Bigr\rangle}

\def\abs#1{\lvert#1\rvert}
\def\Abs#1{\bigl\lvert#1\bigr\rvert}
\def\ABs#1{\Bigl\lvert#1\Bigr\rvert}
\def\ABS#1{\biggl\lvert#1\biggr\rvert}

\def\mid{\,:\,}

 \def\ind#1{ {\mathbf 1}_{#1}}

 \def\sh#1{\operatorname{sh}(#1) }

\def\Xint#1{\mathchoice
   {\XXint\displaystyle\textstyle{#1}}%
   {\XXint\textstyle\scriptstyle{#1}}%
   {\XXint\scriptstyle\scriptscriptstyle{#1}}%
   {\XXint\scriptscriptstyle\scriptscriptstyle{#1}}%
   \!\int}
\def\XXint#1#2#3{{\setbox0=\hbox{$#1{#2#3}{\int}$}
     \vcenter{\hbox{$#2#3$}}\kern-.5\wd0}}

\def\dashint{\Xint-}

 \def\Enl #1,#2,{\operatorname{ Enl}_{#1}(#2)}
 \def\emb #1.#2.{\operatorname{ emb}(#1,#2)}
 \def\zR#1{R_{(#1)}}
 
 \def\dil #1;#2;{\operatorname{Dil}_{#1}#2}
\def\eqdef{\stackrel{\mathrm{def}}{{}={}}}

\begin{document}

\title[Variations on Journ\'e's Lemma]{Variations on the Theme of Journ\'e's Lemma}

\author[C.~Cabrelli]{Carlos Cabrelli}

\email{cabrelli@dm.uba.ar}

\address{Carlos Cabrelli \\
Departamento de Matem\'atica\\
Fac.~de Ciencias Exactas y Naturales, Universidad de Buenos Aires\\
Ciudad Universitaria, Pabell\'on I\\
1428 Capital Federal\\
Argentina}

\author[M.T.~Lacey]{Michael T.~Lacey}

\address{Michael T. Lacey \\
School of Mathematics\\
Georgia Institute of Technology\\
Atlanta,  GA 30332 USA}

\email{lacey@math.gatech.edu}
\thanks{This research is supported in part by an NSF Grant, a Guggenheim Fellowship, and 
the Erwin Schr\"odinger Institute, Vienna Austria. }    
%"

\author[U.~Molter]{Ursula Molter}

\address{Ursula Molter\\
Departamento de Matem\'atica\\
Fac.~de Ciencias Exactas y Naturales, Universidad de Buenos Aires\\
Ciudad Universitaria, Pabell\'on I\\
1428 Capital Federal\\
Argentina}

\email{umolter@dm.uba.ar}
\thanks{The research of the first and third authors is partially supported by grants UBACyT X058, X108, and CONICET PIP 456/98}

\author[J.C.~Pipher]{Jill C.~Pipher }

\address{Jill C.~Pipher\\
Department of Mathematics\\
Brown University\\
Providence RI 02912 USA}

\email{jpipher@math.brown.edu}

\date{}

\begin{abstract}Journ\'e's Lemma \cite{MR87g:42028} is a critical component of many questions related to the 
product $\operatorname{BMO}$ theory of S.-Y.~Chang and R.~Fefferman.  This article presents several
different variants of the 
Lemma, in two and higher parameters, some known, some implicit in the literature, and some new. 
\end{abstract}

\maketitle

%%%%%%%%%%%%%%%%%%%%%%%%%%%%%%%%%%%%%%%%%%%%%%%%%%%%%%%%%%%%%%

\section{Introduction, Journ\'e's Lemma}      \parskip=11pt

We begin the discussion in two dimensions.  Let $\operatorname M$ denote the strong maximal function in the plane. 
Let $\mathcal U$ denote a collection  dyadic rectangles of the plane,  whose union $\sh {\mathcal U}$ is a set
of finite measure, and set $\Enl , \mathcal U,{}\eqdef{}\{\operatorname M\ind {\sh {\mathcal U}} {}>\tfrac12\}$.  
For a dyadic rectangle $R=R_{(1)}\times R_{(2)}\in\mathcal U$, set 
 \begin{equation*} 
\text{emb}(R;\mathcal U){}\eqdef{}\sup\{\mu>1\mid (\mu R_{(1)})\times R_{(2)}\subset \Enl , \mathcal U,\}
 \end{equation*} 
In this display, and throughout this paper, we use the notation $ \lambda R $  to denote the set 
that has the same center as $ R $ but is dilated by an amount $ \lambda $.   (Section~\ref{s.notation} has a comprehensive 
list of notations and conventions.)

The subject of this paper is the result of  J.-L.~Journ\'e  from 1987 \cite{MR87g:42028}. 

 \begin{lemma}\label{l.journe-classic}  [Journ\'e's Lemma]  For all $\epsilon>0$, and any collection $\mathcal U'\subset\mathcal U$ of pairwise incomparable  dyadic rectangles contained in $\mathcal U$, we have the inequality 
 \begin{equation}  \label{e.JOurne}
\sum_{R\in\mathcal U'}\emb R.\mathcal U.^{-\epsilon}\abs {R}\lesssim\abs{\sh {\mathcal U'}}.
 \end{equation}
The implied constant depends only upon $ \epsilon >0 $.
 \end{lemma} 

This Lemma has proven to be an invaluable aid in those problems associated with  
 Carleson measures in product setting.    
In particular, we have been careful to state (\ref{e.JOurne}) as an inequality that is uniform over all 
choices of subset $\mathcal U'\subset\mathcal U$.   We insist on this formulation so that the Lemmas will more readily apply to the 
setting of Carleson measures in the product domain.  The estimation of the norm of subject objects can be quite 
complicated, and the Journ\'e Lemma permits an upper bound in terms of simpler norms. 
See Corollary~\ref{c.BMO}. We comment in more detail on the context of this Lemma in 
the next section. 

 By `product setting' we mean that range of questions which are concerned with 
issues of harmonic analysis that are invariant with respect to a family of dilations with at least two free parameters.

Despite the appearance of this Lemma close to twenty years ago, one cannot yet describe the  precise role that 
this Lemma plays in the product theory, especially when confronting issues related to the induction on the 
number of free parameters.  Indeed, this role 
will be understood by further developments in what seems to be a still nascent product 
theory. Following the work of Chang and Fefferman, see Journ\'e \cite{MR88d:42028}, Carberry 
and Seeger \cite{MR93b:42035}, Fefferman and Pipher \cite{MR1439553},  and the more recent results of Muscalu, Pipher, Tao and Thiele 
\cites{camil1,camil2}, among other papers. We intend this 
paper to be a source book for ideas associated with the Lemma, with a description of what is known, recent 
innovations, as well as some refinements, that as of yet, have not found applications.

There are three themes to the refinements.   First, the Lemma does not appear to 
admit a completely trivial extension to higher parameters.  The point that simplifies the analysis in 
two parameters is that if $R$ and $R'$ are distinct, intersecting rectangles, then it is the case that
two sides of the rectangles are in an inverse relation.  
On the other hand, in three parameters, the different faces of the two rectangles can have a number of 
relations.   See Figure 1 for the situation in the plane.  In fact, the best methods to pass to higher numbers of parameters probably
has not as yet been discovered. 
And there are also versions of the Lemma, for the higher 
parameter setting, in which rectangles are replaced by 
more complicated sets.  These constructions, which are taken up in Section~\ref{s.few} for instance, may in 
applications, permit one to 
get at more directly the particular manner in which e.g.~the $\operatorname{BMO}$ space of three 
parameters is built up from that of two parameters.

	Second, the ``embeddedness term'' $\text{emb}(R,\mathcal U)$ above, is the new element required in two and higher parameters.  There is interest in having different measures of the embeddedness term that are  essentially smaller than that given above. 
	Now, the 	 power of $\epsilon$ in (\ref{e.JOurne}) is used to cancel out  terms that are logarithmic in $\text{emb}(R,\mathcal U)$.  
	Any proof of the Lemma must account for the fact that   a subcollection $\mathcal U'$ in which $\emb R.\mathcal U.\simeq\mu$, for all $R\in\mathcal U'$, one has 
	that $\mathcal U'$ is the union of $O(\log \mu)$ subcollections in which the rectangles are essentially disjoint.

	 If we decrease the embeddedness term, we expect the combinatorial difficulties to multiply, and the logarithmic terms to increase.  We maintain the term $\text{emb}(R,\mathcal U)^{-\epsilon}$, and do not keep track of how  quickly the logarithmic terms increase.    

In typical applications of Journ\'e's Lemma, one obtains, from say ``Schwartz tails'' arguments,  a rapid decrease in 
	 terms of the embeddedness quantity.   

 Third, there are specific instances in which the ``enlarged set'' $\Enl , \mathcal U,$  plays a important role.  As phrased 
 	above, one has $\abs {\Enl ,\mathcal U,}\le{}K\abs {\sh{\mathcal U}} $, with constant $K$ strictly bigger than one. 
	In a paper of Lacey and Ferguson, \cite{MR1961195}, it turns out to be essential that, for arbitrary $\delta>0$, 
	one can select $\Enl , \mathcal U,$ such that $\abs {\Enl , \mathcal U,}\le(1+\delta)\abs {\sh {\mathcal U} } $.  In this regard, 
	also see Lacey and Terwilleger \cite{witherin}.
	We investigate other examples where this can be obtained.

From time to time, we will refer to the set $\Enl , \mathcal U,$ as $V$.  
 It is interesting to note that the conclusion of Journ\'e's Lemma implies the formally stronger conclusion that 
  \begin{equation*}
 \Bigl\lVert \sum_{R\in\mathcal U} \emb R.\mathcal U.^{-\epsilon}\ind R \Bigr\rVert_p\lesssim\abs{\sh {\mathcal U}}^{1/p},\qquad 1<p<\infty.
  \end{equation*}
 This is an immediate consequence of the John Nirenberg inequality, Lemma~\ref{l.jnproduct}, in the product $\operatorname{BMO}$	 setting. 

\smallskip 
%%%%%%%%%%%%%%%%%%%%%%%%%%%%%%%%%%%%%%%%%%%%%%%%%%%%%%%%%%%%%%
\subsection{Notations and Conventions} \label{s.notation}

The dyadic intervals in $\mathbb R$ are 
\begin{equation*}
\mathcal D {}\eqdef{}\{ [j2^k,(j+1)2^k)\mid j,k\in \mathbb Z\}. 
\end{equation*}
This collection of intervals has the \emph{grid property}, namely that for any two intervals $I,I'\in\mathcal D$ it 
is the case that $I\cap I'\in \{\emptyset,I,I'\}$.  We return to the this property below.

  The set $\mathcal D^d$ is then the set of dyadic rectangles 	
	in $d$ dimensional space.  Such a rectangle is a product $R=\prod_{j=1}^d R_{(j)}$, where $R_{(j)}$ 
	is the product in the $j$th coordinate.    $\mathcal U$ denotes a generic subset of $\mathcal D^d$.  The {\em shadow} of 
	$\mathcal U$ is 
	 \begin{equation*}
	\sh {\mathcal U}{}\eqdef{}\bigcup_{R\in\mathcal U}R.
	 \end{equation*}

For a rectangle $R$ and $\lambda >0$ we set $\lambda R$ to be the rectangle with the same center as $R$, and 	
	whose dimension in each coordinate are to be $\lambda$ times the corresponding dimension of $R$. 
It will be useful to have a more versatile notion of dilations.  Thus,  for a vector $\vec \lambda=(\lambda_1,\ldots,\lambda_d)\in\mathbb R^d$, 
set 
 \begin{equation} \label{e.dil} 
\dil \vec \lambda;R;{}\eqdef{}\otimes_{j=1}^d \lambda_j R_{(j)}
 \end{equation}

The dyadic rectangles are ordered by inclusion;  {\em maximal  elements } of $\mathcal U$ refer to rectangles 
	that are maximal with respect to inclusion.   This is quite a good partial order in one dimension of 
	course: Two dyadic intervals intersect iff they are related under the partial order.   
	It is less effective  in higher dimensions, though a distinguishing feature of two dimensions is that 
	if two non equal rectangles intersect and are not comparable,
	then the two sides of the rectangles must be in reverse order with 
	respect to inclusion.      This fails  three parameters, and explains in part the difficulty 
	in moving from two to three parameters in some of our arguments.

A set of dyadic rectangle $\mathcal U$ has {\em scales separated by $\mu$} iff for any two rectangles 
	$R,R'\in\mathcal U$, and for any $j$, if $\abs{R_{(j)}}<\abs{R'_{(j)}}$ then $\mu\abs{R_{(j)}}<\abs{R'_{(j)}}$.   Any set $\mathcal U$ is a union of $\lesssim(\log \mu)^d$ subsets 
	which have scales separated by $\mu$. 
	An example fact we shall rely upon is this.  If $I$ and $J$ are   intervals, with   $\abs I<{}\abs J$ and 
	$I\cap J\not=\emptyset$,  then it is the case that 
	 \begin{equation} \label{e.expanding}
	 I\subset(1+\tfrac {\abs I}{\abs J})J\subset\{ \operatorname  M\ind J>(1+2\tfrac {\abs I}{\abs J})^{-1}\}
	 \end{equation}
	where $\operatorname M$ is the maximal function in one dimension. We will be applying this with $I$ and $J$ dyadic intervals, 
	and with scales separated by some large amount.

The strong maximal function is 
	 \begin{equation*}
	\operatorname M f(x)=\sup_{x\ni R}\dashint_R \abs {f(y)}\; dy,
	 \end{equation*}
	where the supremum is taken over all (non--dyadic) rectangles $R$ in $\mathbb R^d$.  In addition, we  use the notation 
	 \begin{equation*}
	\dashint_A f\; dx=\abs{A}^{-1}\int_A f\; dx.
	 \end{equation*}
	We use without comment the $L^p$ inequalities known for the strong maximal function.
	For intersecting rectangles $R$ and $R'$ we have 
	 \begin{equation*}
	R\subset \prod_{j=1}^d  \gamma_j \zR j'\subset \bigl\{ \operatorname M  \ind {R'}>\prod_{j=1}^d \gamma ^{-1}_j\Bigr\}
	 \end{equation*}
	where $\gamma _j{}\eqdef{}1+\abs{\zR j}\abs{\zR j'}^{-1}$.

	It is known, see e.g.~the work of Melas\footnote{In fact, Melas computes the 
	exact constant in the weak type inequality.} \cites{MR2003c:42021,MR1973058}, 
	that even in one dimension,  the maximal function maps $L^1$ into $L^{1,\infty}$ with norm strictly bigger than one.  
 The dyadic maximal function however maps $L^1$ into $L^{1,\infty}$ with norm $1$. 
	 We shall have need of a  variant of  this
	 well--known fact.  
	 
	 Define a {\em grid} to be a collection
	 $\mathcal I$  of intervals in the real line  for which for all $I,I'\in\mathcal I$, 
	 $I\cap I'\in\{\emptyset,I,I'\}$.  For a collection of intervals $\mathcal I$, not necessarily
	 a grid,  set 
	  \begin{equation} \label{e.M^I}
	\operatorname M^\mathcal I f(x){}\eqdef{}\sup_{I\in\mathcal I} \ind I(x)\dashint_I f(y)\; dy.
	  \end{equation}
	Then, for any grid $\mathcal I$,  $\operatorname M^\mathcal I$ maps $L^1(\mathbb R)$ into into $L^{1,\infty}(\mathbb R)$ with norm one.  
	This, in
	particular, is  true for the dyadic grid $\mathcal D$.   
	
	This fact we prove here, for the sake of completeness.  For a non negative integrable $f$, and $\lambda>0$, 
	the set $\{\operatorname M^\mathcal I f>\lambda\}$ is union of intervals in the grid.  Hence is a disjoint union of 
	intervals in $\mathcal I'\subset \mathcal I$.  For each interval $I\in\mathcal I'$, we must have 
	 \begin{equation*}
	\dashint_I f\; dx\ge\lambda .
	 \end{equation*}
	Hence, 
	 \begin{equation*}
	\lambda \abs{\{\operatorname M^\mathcal I f>\lambda\}}=\lambda\sum_{I\in\mathcal I'}\abs {I}\le{}\sum_{I\in\mathcal I'}\int_I f \; dx\le{}\norm f.1..
	 \end{equation*}

We shall have need of a notion of {\em shifted dyadic grids},  due to M.~Christ,  defined as follows. 
The definition of the grids depends upon a choice of integer $\mathsf d$, and set  $\delta=(2^{\mathsf d}+1)^{-1}$ for integer $\mathsf d$.   For integers
$0\le{}b<{}\mathsf d$, and $\alpha\in\{\pm(2^{\mathsf d}+1)^{-1}\}$, let 
  \begin{equation} \label{e.shifted-grids}
\begin{split}
 \mathcal D_{\mathsf d,b,\alpha}&{}{}\eqdef{}\{2^{k\mathsf d+b}((0,1)+j+(-1)^k\alpha)\mid k\in\mathbb Z,\ j\in\mathbb Z\}.
 \\
 \mathcal D_{\mathsf d}&{}{}\eqdef{}\bigcup_\alpha\bigcup_{b=0}^{\mathsf d-1}\mathcal D_{\mathsf d,b,\alpha}.
 \end{split}
  \end{equation}
 One checks that $\mathcal D_{\mathsf d,b,\alpha}$  is a grid.  Indeed, it suffices to assume
 $\alpha=(2^{\mathsf  d}+1)^{-1}$, and that $b=0$.  Checking the grid structure can be done by
 induction.  And it suffices to check that the intervals in $\mathcal D_{\mathsf d,0,\alpha}$ of
 length one are a union of intervals in $\mathcal D_{\mathsf d,0,\alpha}$ of length $2^{-\mathsf d}$. 
 One need only check this for the interval $(0,1)+\alpha$.  But certainly
  \begin{align*}
 (0,1)+{(2^{\mathsf d}+1)}
 {}&{}=\bigcup_{j=0}^{2^{\mathsf d}-1} (0,2^{-d})+\frac j{2^{\mathsf d}}+{(2^{\mathsf d}+1)}
 \\&{}={}
 \bigcup_{j=0}^{2^{\mathsf d}-1} (0,2^{-d})+\frac {j+1}{2^{\mathsf d}}+{2^{\mathsf d}(2^{\mathsf d}+1)}
  \end{align*}
 And this proves the claim.

   What is just as  important concerns the collections
 $\mathcal D_{\mathsf d}$.  For each dyadic interval
 $I\in\mathcal D$, $I\pm\delta\abs I\in\mathcal D_{\mathsf d}$.\footnote{The problem we are avoiding here is that the dyadic grid distinguishes  dyadic rational points. 
 At the point $0$, for instance,
 observe   that for all integers $k$, $(1+\delta)(0,1)\not\subset(0,2^k)$,
 regardless of how big $k$ is.}   Moreover, the maximal function $\operatorname M^{\mathcal D_{\mathsf d}}$ maps $L^1$ into
 $L^{1,\infty}$\footnote{In fact, taking $\mathsf d=1$, it is routine to check that $\operatorname M^{\mathcal D_{\mathsf 1}}$ dominates 
an absolute multiple of the usual maximal function.  Thus, proving that it satisfies the weak type inequality.} 
 with norm at most $2\mathsf d\simeq\abs{\log \delta}$. In fact we need the finer estimate
 \begin{equation}\label{e.small-weak} 
\abs{ \{\operatorname M^{\mathcal D_{\mathsf d}}\ind {\sh{U}}>1-\delta\}}\le{}(1+K\delta\mathsf d)\abs{\sh{U}},
 \end{equation}
for all subset $U$ of the real line of finite measure, and some constant $K$.  
This is an effective estimate since $\delta\mathsf d\simeq \delta\abs{\log \delta}\longrightarrow0$, 
as $\delta\longrightarrow0$.  

To see this estimate,  note that 
 \begin{align*}
\abs{\{ \operatorname M^{\mathcal D_{\mathsf d}}\ind {\sh{U}} >1-\delta\}}\le{}&
\abs{\sh{U}}+\sum_{b=0}^{d-1}\sum_{\alpha\in \{\pm (2^{\mathsf d}+1)^{-1}\}} 
\abs{ U^c\cap \{ \operatorname M^{\mathcal D_{\mathsf d, b, \alpha}}\ind {\sh{U}} >1-\delta\}}
\\{}\le{}&(1+2\mathsf d[(1-\delta)^{-1}-1])\abs{\sh{U}}. 
 \end{align*}

In statements of Journ\'e's Lemma, $\mathcal U$ will denote a generic collection of rectangles  of $\mathbb R^d$, whose shadow is  of finite measure. 
    The statement of the Lemma will depend upon a particular choice of  {\em enlarged set}
    which we will always define in terms of some maximal function.  It will be denoted as 
    $\Enl ,\mathcal U,$, or more simply as $V$.  At times, this definition will  be iterated.  In this case, we denote the enlarged set as   $\Enl j,U,$, the subscript 
    $j$ denoting     the number of times the definition is iterated.\footnote{In some statements of Journ\'e's Lemma, the role of the enlarged 
    set is suppressed, and only the ``embeddedness'' terms are used. In this paper, we are of course concerned with the selection of the ``enlargement'' and some of the enlargement's properties.}  

Journ\'e's Lemma also depends upon a notion of {\em embeddedness} of a rectangle $R\in\mathcal U$, relative 
	to the enlarged set $\Enl ,\mathcal U,$.  If there is no ambiguity about the enlarged set, we use the notation 
	$\emb R.\mathcal U.$.  Otherwise, the notation $\emb R. {\Enl j,U,}.$ is used.   The definition of the enlarged set, and the notion of embeddedness will vary from section from section, but the notation will not. 

Many factors, arising in most instances from combinatorial considerations,  are increasing like 
	a power of $\log \emb R.\mathcal U.$.  These terms are considered to be inconsequential.   One instance of this
	which frequently arises is as follows.    Let $\mu>1$ and let $\mathcal U$ be a collection of rectangles  with a 
	shadow of finite measure.  Let $\mathcal U'\subset\mathcal U$ satisfy  $\mu\le \emb
	R.\mathcal U.\le2\mu$ for all $R\in\mathcal U'$, and 
	the scales of $\mathcal U'$ are separated by $10^{3d}\mu$.  Then to prove Journ\'e's Lemma, it suffices to show that 
	 \begin{equation} \label{e.standard}
	\sum_{R\in\mathcal U'}\abs R\lesssim\abs {\sh {\mathcal U'}}.
	 \end{equation}
	This we will refer to as the {\em standard reduction.}  
	
	This last inequality obviously holds if  the rectangles in $\mathcal U$ 
	are {\em essentially disjoint.}  That is, there is a choice of absolute constant $c$, and there 
	are sets $E(R)\subset R$ so that 
	$\{E(R)\mid R\in\mathcal U\}$ are pairwise disjoint sets.  And that $\abs {E(R)}\ge{}c\abs R$.  
 Obtaining this, or a property similar to it, 
	is an obvious strategy for proving (\ref{e.standard}) in a manner that is uniform with respect
	to $\mathcal U'\subset\mathcal U$.

We write $A\lesssim{}B$ if there is an (unimportant)  absolute constant $K$ (permitted to depend upon parameters 
as specified in e.g.~the statement of a Proposition) such that $A\le{}KB$.   $A\simeq B$ means that
$A\lesssim{}B$ and $B\lesssim{}A$.

%%%%%%%%%%%%%%%%%%%%%%%%%%%%%%%%%%%%%%%%%%%%%%%%%%%%%%%%%%%%%%

\section{Hardy Space, $\text{BMO}$, and Carleson Measures in the Product Theory}

The realm of application of Journ\'e's Lemma is to the product $\operatorname{BMO}$ theory. 
Functions in this class are described 
by their Carleson measures.   We survey these subjects, beginning with the Carleson measures, and including 
explicit definitions and Lemmas that follow from the more purely geometric versions of Journ\'e's
Lemma that are in other parts of this paper.

\subsection{Carleson Measures}
Journ\'e's Lemma is most directly applied to the control of Carleson measures in the product setting.  
And we first address this implication, following up with connections to the product Hardy space theory.  

For a map $\alpha\mid \mathcal D^d\longrightarrow\mathbb R_+$, set 
 \begin{equation} \label{e.CM} 
\norm \alpha .CM. {}\eqdef{}\sup_{\mathcal U}  \abs{\sh{\mathcal U}}^{-1} \sum_{R\in\mathcal U}{\alpha(R)}  .
 \end{equation}
``CM'' is for Carleson measure. What is most essential here is that the supremum is taken over all 
subsets $U\subset\mathbb R^d$ of finite measure.     In one dimension, a small additional argument permits one to  restrict the supremum to intervals.   

This definition is confusing, as there are no measures present.  In  Section~\ref{s.cm} 
we recall the more classical definition of a Carleson measure.

In more parameters, it is natural to suppose  that one should be able to restrict the supremum above to rectangles.  While this is not the case,\footnote{Historically, these examples did 
not arrive in this way, but where phrased in the language of the Hardy space $\operatorname H^1$, and it's dual. 
We comment in more detail below.}
this supremum does play a distinguished role in the theory, and we denote this supremum by $\norm \alpha. CM(\text{rec}).$.  

In particular, in dimensions $2$ and higher, for all $\epsilon>0$,  there are Carleson measures $\alpha$ with $\norm \alpha.CM.=1$ and $\norm \alpha.CM(\text{rec}).<\epsilon$.  
The main application of  Journ\'e's Lemma is to show that despite this general difficulty, we can in some instances use the rectangular norm to control the general norm.  

 \begin{corollary}\label{c.CM}  For all $\epsilon>0$,  all $\mu>1$, and collections of rectangles $\mathcal U$ whose shadow has finite area in the plane,  let $\mathcal U_\mu\subset\mathcal U$ be a collection of rectangles with $\text{\rm emb}(R, \mathcal U)\simeq\mu$.  Then, 
 \begin{equation*}
\norm \alpha|_{\mathcal U_\mu}.CM.\lesssim\mu^{\epsilon}\norm \alpha .CM(\text{\rm rec}).
 \end{equation*}
 \end{corollary} 

It is to be stressed that  this Lemma, as stated, is restricted to the plane.  With more than two parameters, we need to either 
take more care with the definition of embeddedness, or with the definition of the ``rectangular'' norm.

\begin{proof}  We should see that for all sets $\mathcal V\subset\mathcal U_\mu$, we have 
 \begin{equation*}
\sum_{R\in\mathcal V}\alpha(R)\lesssim{}\mu^{\epsilon}\norm \alpha .CM(\text{\rm rec}).\abs {\sh {\mathcal V}} .
 \end{equation*}
Let $\mathcal V'$ be the maximal dyadic rectangles in $\mathcal V$.  Then, 
 \begin{align*}
\sum_{R'\in\mathcal V'} \sum_{\substack{R\in\mathcal V\\ R\subset R'}} \alpha(R)\le{}& 
	\norm \alpha .CM(\text{\rm rec}).\sum_{R'\in\mathcal V'}\abs {R'}
	\\{}\lesssim{}&\mu^{\epsilon}\norm \alpha .CM(\text{\rm rec}).\abs {\sh {\mathcal V}} 
	 \end{align*}
In the top line we have used the definition of the rectangular Carleson measure norm, and in the bottom 
Journ\'e's Lemma, as stated in Lemma~\ref{l.journe-classic} say.
\end{proof}

 In higher parameters, one can continue to use the rectangular norm, using instead of the planar version of Journ\'e's Lemma, the 
 form as stated in Lemma~\ref{l.journe-rectangles}.  What is more interesting is to define a notion of Carleson measure norms that 
 uses  Lemma~\ref{l.few}.  Towards this end, let us say that a collection $\mathcal U$ of rectangles in $\mathbb R^d$ has $\ell$ parameters
 iff there is a subset $L\subset\{1,\ldots,d\}$ with $\abs L=\ell$, so that for any two rectangles $R,R'\in\mathcal U$ we have $
 \zR j=\zR j'$ for all $j\not\in L$.  
 Let us set 
  \begin{equation*}
 \norm \alpha .CM(\ell).=\sup_{\text{$\mathcal U$ $\ell$ parameters}} 
 \abs{\sh {\mathcal U} }^{-1} \sum_{R\in\mathcal U} \alpha (R). 
  \end{equation*}
 Notice that the $CM(1)$ norm reduces to essentially the most natural extension of the 
 rectangular norm to higher parameters.    These norms will increase in $\ell$. In 
 general, one cannot control the $CM(\ell)$ norm by the $CM(\ell-1)$ norm, except through devices 
 like Journ\'e's Lemma.
 
 This definition goes someway towards capturing the subtle way that Carleson measures of $d$ parameters are built up from 
 those of $d-1$ parameters.  In particular, we have the following Lemma, in which we use the notations of (\ref{e.d-11}) and (\ref{e.d-12}). 
 We only state this Lemma in the case of $\ell=d-1$ as it is the only case that has found application to date.

  \begin{proposition}\label{p.cm-d-1}  For all $\delta>0$ the following holds. 
 Let $\mathcal U$ be a collection of rectangles in $\mathbb R^d$ whose shadow has finite measure, and for $\mu>1$ set 
  \begin{equation*}
 \mathcal U_\mu{}\eqdef{}\{R\in\mathcal U\mid \mu\le{} \emb R. \sh {\mathcal U} .\le{}2\mu \}.
  \end{equation*}
 Then, we have 
  \begin{equation*}
 \norm \alpha| _{\mathcal U_\mu} .CM(d).\lesssim{}\mu^\epsilon\norm \alpha .CM(d-1).
  \end{equation*}
 The implied constant depends upon   $\epsilon>0$. 
  \end{proposition} 
 
 These concepts, and this lemma are used in Lacey and Terwilleger \cite{witherin}.

\bigskip 
 An important aspect of the subject is  the connection of the definition of 
 Carleson measures to a John--Nirenberg inequality.

 \begin{lemma}\label{l.jnproduct}  We have the inequality below, valid for all  collections of rectangles $\mathcal U$ whose shadows have finite measure.
 \begin{equation*}
\NOrm \sum_{R\in\mathcal U} \frac {\alpha(R)}{\abs R}\ind R .p.\lesssim{} \norm \alpha .CM. \abs{\sh{\mathcal U}}^{1/p},\qquad 1<p<\infty.
 \end{equation*}
 \end{lemma}

\begin{proof}  This is the proof by duality from \cite{MR82a:32009}. Let $\norm \alpha.CM.=1$. 
Define 
 \begin{equation*}
F_V{}\eqdef{}{}
\sum_{R\subset V} \frac {\alpha(R)}{\abs R}\ind R 
 \end{equation*}
We shall show that for all $\mathcal U$, there is a set $V$ satisfying $\abs{V}<2\abs{\sh{\mathcal U}}$ for which 
 \begin{equation} \label{e.UV}
\norm F_{\sh{\mathcal U}} .p.\lesssim{} \abs{\sh{\mathcal U}}^{1/p}+\norm F_V .p.
 \end{equation}
Clearly, inductive application of this inequality will prove our Lemma.

The argument for (\ref{e.UV}) is by duality.  Thus, for a given $1<p<\infty$, and conjugate index $p'$, take $g\in
L^{p'}$ of norm one so that $\norm F_U.p.=\ip F_U,g.$.  Set 
 \begin{equation*}
V=\{ \operatorname Mg >K \abs{ \sh{\mathcal U}}^{-1/p'}\}
 \end{equation*}
where $\operatorname M$ is the strong maximal function and $K$ is sufficiently large so that $\abs V<2\abs{\sh{\mathcal U}}$.  Then, 
 \begin{equation*}
\ip F_{\sh{\mathcal U}},g.=\sum_{\substack{R\in\mathcal U\\R\not\subset V}} \alpha(R)\dashint_R g \; dx+\ip F_V ,g.
 \end{equation*}
The second term is at most $\norm F_V.p.$ by H\"older's inequality. 
%"
For the first term, note that the 
average of $g$ over $R$ can be at most $K\abs{\sh{\mathcal U}}^{-1/p'}$.  So by the definition of Carleson measure norm, it is at most 
 \begin{equation*}
\sum_{\substack{R\in\mathcal U\\R\not\subset V}} \alpha(R)\dashint_R g \; dx\lesssim{}\abs{\sh{\mathcal U}}^{-1/p'}\sum_{R\in\mathcal U}\alpha(R)\lesssim{}\abs{\sh{\mathcal U}}^{1/p},
 \end{equation*}
as required by (\ref{e.UV}). 

\end{proof}

%%%%%%%%%%%%%%%%%%%%%%%%%%%%%%%%%%%%%%%%
\subsection{Classical Definition, Carleson Embedding Theorem}   \label{s.cm}

Our use of the the term ``Carleson measure'' is not the standard one. Given a function $\alpha\mid \mathcal D^d\longrightarrow\mathbb R_+$, define a 
measure on $\mathbb R^d\times \mathbb R_+^d$ by 
 \begin{equation*}
\mu_\alpha=\sum_{R\in\mathcal D^d} \alpha(R) \delta_{ R\times\norm R..}
 \end{equation*}
where $\norm R..=(\abs{R_{(1)}},\ldots,\abs{R_{(d)}})$.  In the instance that $\alpha(R)=\abs{R}^{-1}\abs{\ip f,h_R.}^2$, the measure 
$\mu_\alpha$ is of the type associated with the area integral of $f$.  (Indeed, in this setting both the continuous and discrete 
formulations are equivalent.)

 For a set $U\subset \mathbb R^d$, define an associated set $\operatorname{Tent}(U)\subset \mathbb R^d\times \mathbb R_+^d$ by 
  \begin{equation*}
 \operatorname{Tent}(U){}\eqdef{}\bigcup_{\substack{ R\in\mathcal D^d \\ R\subset U} } R \times [0,\abs{R_{(1)}}]\times\cdots\times[0,\abs{R_{(d)}}]
  \end{equation*}
 This is the {\em tent over $U$.}
 Then, the substance of the Carleson measure condition is the inequality 
  \begin{equation*}
 \mu_\alpha(\operatorname{Tent}(U))\le{} \norm \alpha .CM. \abs{U} ,
  \end{equation*}
 for all sets $U\subset\mathbb R^d$ of finite measure.  
 Notice that the left hand side concerns objects of $2d$ dimensions, while the right hand side has only dimension $d$.

 The importance of the Carleson measure condition arises from the Carleson Embedding Theorem, which we again state 
 in a discrete form.  Given a function $\alpha\mid \mathcal D^d\mapsto [0,\infty) $, define an operator 
  \begin{equation*}
 \operatorname T _{\alpha} f{}\eqdef{}\sum _{I\in\mathcal D^d} \alpha(R) \mathbf 1 _{R }\dashint _R f(y)\; dy 
  \end{equation*}

 \begin{theorem}\label{t.embed}   We have the equivalence below, valid for all $1<p<\infty $.
 \begin{equation*}
\norm \operatorname T_\alpha.p.\simeq{} \norm \alpha.CM. 
 \end{equation*}
 \end{theorem}

\begin{proof}
The inequality $\norm T_\alpha.p. {}\gtrsim{}{} \norm \alpha.CM. $ follows by testing the operator  $\operatorname T_\alpha $
against a function $f=\mathbf 1 _{U} $.  Thus, 
 \begin{align*}
\NOrm \sum _{R\subset U} \alpha(R) \mathbf 1 _{R} .p.&{}\le{}\norm \operatorname T _{\alpha} \mathbf 1 _{U}.p.
\\
&{} {}\lesssim{}\norm \operatorname T_\alpha.p.\abs U ^{1/p}.
 \end{align*}
This condition appears stronger than the definition of the Carleson measure norm, but the John--Nirenberg inequality of 
course implies that it is equivalent to this definition.

And we shall find  the John Nirenberg inequality essential for the proof of the reverse inequality. 
We do not prove the strong type inequality directly, but rather prove the weak type inequality 
 \begin{equation*}
\abs{ \{ \operatorname T _\alpha f>\lambda\} } {}\lesssim{}\norm \alpha.CM.^p \lambda^{-p} \norm f.p.^p,\qquad 1<p<\infty.
 \end{equation*}
To prove this inequality, 
let us observe that the definition of the Carleson measure norm, and that of the operators $T _{\alpha } $ is invariant 
under dilations.  Namely, letting $\mu $ be  any power of $2 $, and setting 
 \begin{equation*}
\beta(R){}\eqdef{}\alpha(\dil (\mu,\cdots,\mu); R;) 
 \end{equation*}
we have $\norm \beta.CM.=\norm \alpha.CM. $.  And,  
 \begin{equation*}
\operatorname T _\beta= \operatorname T _\alpha \dil (\mu,\cdots,\mu); 1;
 \end{equation*}
Thus, it suffices to prove a single instance of the weak type inequality.  Namely, 
that for $1<p<\infty $, there is a constant $K_p $ so that for all $\alpha $ with Carleson measure norm $1 $, and all 
functions $f\in L^p $ of norm one, we have 
 \begin{equation}  \label{e.weak-instance}
\abs{ \{ \operatorname T _{\alpha } f>1\} }\le{}K_p.
 \end{equation}

We inductively decompose the collection of dyadic rectangles. In the base step, take 
 \begin{equation*}
\mathcal U_0{}\eqdef{}\{ R\in\mathcal D^d\mid  \dashint_R f(y)\; dy\ge1\} 
 \end{equation*}
Set $\operatorname {Stock}{}\eqdef{}\mathcal D^d-\mathcal U_0 $.  In the inductive stage, given $\mathcal U_0,\cdots,\mathcal U_k $, 
 to construct $\mathcal U_{k+1} $, we set 
  \begin{equation*}
 \mathcal U_{k+1}{}\eqdef{}\{ R\in \operatorname {Stock}\mid  \dashint_R f(y)\; dy\ge2 ^{-k+1}\}. 
  \end{equation*}
Then, update $\operatorname {Stock}{}\eqdef{}\operatorname {Stock}-\mathcal U _{k+1} $. 

By the strong Maximal Function estimate, we have 
 \begin{equation*}
\abs{ \sh {\mathcal U_k} } {}\lesssim{}2^{kp},\qquad k\ge0. 
 \end{equation*}
Thus, we shall not even estimate $\operatorname T _\alpha f $ on the set $\sh {\mathcal U_0} $.  

For the collections $\mathcal U_k $ for $k\ge1 $, we have an upper bound on the average of $f $ over those 
rectangles $R\in\mathcal U_k $.   This, with the John--Nirenberg, will give us a 
favorable estimate in $L^s $ norm, for a choice of $s>p $.  
 \begin{align*}
\NOrm \sum _{R\in\mathcal U_k} \alpha(R) \dashint _R f (y)\; dy \mathbf 1 _{R}.s. 
	&\le{}2^{-k+1}\NOrm \sum _{R\in\mathcal U_k} \alpha(R)  \mathbf 1 _{R}.s.
\\
	& {}\lesssim{}2^{-k} \abs{\sh {\mathcal U_k} } ^{1/s}
\\
	& {}\lesssim{}2 ^{-k(1-p/s)}.
	 \end{align*}
This is summable over $k\ge1 $, and so easily completes the proof of (\ref{e.weak-instance}).

\end{proof}

%%%%%%%%%%%%%%%%%%%%%%%%%%%%%%%%%%%%%%%%%%%%%%%%%%%%%%
\subsection{The Product Hardy Theory} \label{s.hardy} \parskip=11pt
We turn to the product Hardy space theory, as developed by S.-Y.~Chang and R.~Fefferman \cites{MR86g:42038
,
MR82a:32009
,
MR90e:42030
,
MR86f:32004
,
MR81c:32016}.   $\operatorname H^1(\mathbb C_+^d)$ will denote  the real  $d$--fold product Hardy space.  This space consists of functions $f\mid 
\mathbb R^d\longrightarrow\mathbb R$.  $\mathbb R^d$ is viewed as the boundary of 
 \begin{equation*}
\mathbb C_+^d=\prod_{j=1}^d \{z\in\mathbb C\mid \operatorname{Re}(z)>0\}
 \end{equation*}
And we require that there is a function $F\,:\,\mathbb C_+^d\longrightarrow\mathbb C$ that is holomorphic in each variable separately, 
and 
 \begin{equation*}
f(x)=\lim_{\norm y..\to0}\operatorname {Re}(F(x_1+iy_1,\ldots,x_d+iy_d)).
 \end{equation*}	
The norm of $f$ is taken to be 
 \begin{equation*}
\norm f. \operatorname H^1 .=\lim_{y_1\downarrow 0}\cdots\lim_{y_d\downarrow 0}\norm F(x_1+y_1,\ldots,x_d+y_d).L^1(\mathbb R^d).
 \end{equation*}

The dual of this space is $\operatorname {Re}\operatorname H^1(\mathbb C_+^d)^*=\operatorname{BMO}(\mathbb C_+^d)$, the $d$--fold product $\operatorname{BMO}$ space. It is a Theorem of S.-Y.~Chang and R.~Fefferman 
\cite{MR82a:32009} that this space has a characterization in terms of the product Carleson measure introduced above. 
We need the product Haar basis.  Thus, set 
 \begin{equation*}
h(x)=-\ind {[-\tfrac12,0]}(x)+\ind {[0,\tfrac12]}(x),\qquad h_I(x)=h\biggl( \frac{x-c(I)}{\abs I}\biggr),\qquad I\in\mathcal D.
 \end{equation*}
The functions $\{h_I\mid I\in\mathcal D\}$ are the Haar basis for $L^2(\mathbb R)$, which is closely associated with the analysis of singular integrals.
 For a rectangle $R=\prod_{j=1}^d R_{(j)}\in\mathcal D^d$
set 
 \begin{equation*}
h_R(x_1,\ldots,x_d)=\prod_{j=1}^d h_{R_{(j)}}(x_j).
 \end{equation*}
The basis $\{ h_R\mid R\in\mathcal D^d\}$ is the $d$--fold tensor product of the Haar basis.  Then it is the Theorem of Chang and Fefferman that
 the product $\operatorname{BMO}$ space has the equivalent norm 
 \begin{align*}
\lVert b\rVert_{\operatorname{BMO}}^2={}&\sup_{\mathcal U\subset\mathbb R^d} \abs{\sh{\mathcal U}}^{-1} \sum_{R\in\mathcal U}\abs{\ip b,h_R.}^2
\\
{}={}& \norm R\longrightarrow \abs{\ip b,h_R.}^2 .CM.
 \end{align*}

The next comments are specific to the two parameter case, namely 
$\operatorname H^1(\mathbb C_+^2)$. 
The space $\operatorname{BMO}(\text{rec})$ has the definition 
 \begin{equation*}
\lVert b\rVert_{\operatorname{BMO}(\text{rec})}^2={} \norm R\longrightarrow \abs{\ip b,h_R.}^2 .CM(\text{rec}).
 \end{equation*}
It was at first, natural supposition that this space is the dual to $\operatorname H^1$. This stems in part from the fact that the rectangular 
$\operatorname{BMO}$ norm has an equivalent formulation in terms that look quite familiar: 
 \begin{gather*}
\norm b. \operatorname{BMO}(\text{rec}).^2={}
\\
\sup_{I\times J}\dashint_{I\times J}  
\ABs{b (x,y)-\dashint_I f(x,y)\;dx -\dashint_J f(x,y)\; dy+\dashint_{I\times J} f(x,y)\; dx,dx}^2\; dx\,dy.
 \end{gather*}
This of course looks like the familiar intrinsic definition of $\operatorname{BMO}$ in terms of bounded mean oscillation over 
intervals in the real line. 
But an example of Carleson  \cite{carleson-example} 
consisted of a class of functions which acted as linear functionals on $\operatorname H^1$ with norm one, yet 
had arbitrarily small $\operatorname{BMO}(\text{rec})$ norm.  This example is recounted at the beginning of 
R.~Fefferman's article \cite{MR81c:32016}.

Parallel to Corollary~\ref{c.CM}, we have this Corollary to Journ\'e's Lemma. 

 \begin{corollary}\label{c.BMO}   For all $\epsilon>0$,  all $\mu>1$,  collections $\mathcal U$ of rectangles in the plane whose shadow has finite measure, let $\mathcal U_\mu$ be a collection of rectangles with $\text{\rm emb}(R; U)\simeq\mu$.  Then, 
 \begin{equation*}
\NOrm \sum_{R\in\mathcal U_\mu}\langle f, h_R\rangle h_R .\operatorname{BMO}.\lesssim\mu^{\epsilon}\norm b .\operatorname{BMO}(\text{\rm rec}).
 \end{equation*}
 \end{corollary}

There is a corresponding notion of a $\operatorname{BMO}(d-1)$ norm, and a Lemma that is parallel to Proposition~\ref{p.cm-d-1}, but we will not state it explicitly.

\medskip

%%%%%%%%%%%%%%%%%%%%%%%%%%%%%%%%%%%%%%%%%%%%%%%%%%%%%%%%%%%%%%
\section{Journ\'e's Lemma in Two Parameters}\parskip=11pt  

We state and prove different versions  of Journ\'e's Lemma in the two parameter setting.
  We 	shall be explicit about the definition of the expanded set, and somewhat flagrant with logarithms of 
	$\emb R.\mathcal U.$.  This is in contrast to the original references, which give slightly more precise 
	estimates for the sum in the Lemma than we do.  

\subsection{The Original Formulation} 

Two proofs of the Lemma in its original formulation, namely Lemma~\ref{l.journe-classic}, are given. 

\subsubsection{The First Proof.} 
We define, as above, $\Enl , \mathcal U,{}\eqdef{}\{ \operatorname M\ind {\sh{\mathcal U}} >\tfrac12\}$, and 
 \begin{equation} \label{e.1st-emb}
\emb R.\mathcal U.{}\eqdef{}\sup\{ \mu\ge1\mid \dil (\mu,1);R;\subset \Enl , \mathcal U,\}
 \end{equation}
In particular we only expand $R$ in it's first coordinate.  We are to prove  Lemma~\ref{l.journe-classic}. 

We pass to the standard reduction,\footnote{It will be clear that in this instance we need only separate 
	scales in first coordinate, not both as we have defined the standard reduction. }
	see (\ref{e.standard}).   In particular we will use the ``essentially disjoint'' argument mentioned immediately 
	below (\ref{e.standard}).
	Say that $R<_1R'$ if $R\cap R'\not=\emptyset$ and $\abs{ R_{(1)}}<\abs {R_{(1)}'}$.  
	Consider the collection $\text{Bad}$ of rectangles $R\in\mathcal U$ for which there are $R^1,\ldots, R^J$ 
	in $\mathcal U$ with $R<_1R^j$ and finally, that 
	 \begin{equation*}
	\ABs{ R\cap \bigcup_{j=1}^J R^j }\ge\tfrac78\abs R .
	 \end{equation*}

	Observe that the collection $\text{Bad}$ is empty.  Indeed, if $R\in\text{Bad}$, then it must be the case that $\emb R.\mathcal U.\ge10\mu$, which is a contradiction.  This   is a straightforward consequence of (\ref{e.expanding}) and the fact that we defined the enlarged set in terms of the  
	strong maximal function. 
	Thus, there is at least $\frac18$ of each rectangle $R\in\mathcal U$ that is disjoint from all other rectangles in $\mathcal U$. 
And so the rectangles in $\mathcal U$ are essentially disjoint, and we have completed the proof. 

\subsubsection{The Second Proof.} 
The second proof  begins with a key new definition  for dyadic intervals $I $ and integers $k\ge0 $. 
For a subcollection  $\mathcal U'\subset\mathcal U $ that is fixed, set 
 \begin{equation*}
\mathcal E(I,k){}\eqdef{}{}\{ I\times J\in\mathcal U' \mid \operatorname {Dil} _{(2^k,1)} I\times J \subset \sh {\mathcal U} \}.
 \end{equation*}
We will suppress the dependence on the choice of $\mathcal U' $. 
There are two points to observe.  First, due to the maximality of the rectangles in $\mathcal U $, we have 
 \begin{equation} \label{e.2nd} 
\sum _{R\in\mathcal E(I,k) }\abs R {}\lesssim{}2\abs{\sh{\mathcal E(I,k) } }. 
 \end{equation}
Second, we consider two dyadic intervals $I\subset I' $, with $2^n\abs I\le\abs {I'} $.  If it is the  case 
that  
 \begin{equation} \label{e.22nd} 
\abs{I\times J\cap \sh{\mathcal E(I',k) } }\ge\tfrac12\abs {I\times J }, 
 \end{equation}
then it must be that $\emb I\times J.\mathcal U .\ge2^{k+n-1} $.

Thus, if we take $\mathcal U' $ to be a subset of rectangles $R\in\mathcal U $, with $\emb R .\mathcal U .\le2^{k_0}  $ for all rectangles.  
It follows from (\ref{e.22nd}) that 
 \begin{equation*}
\sum _{I \in\mathcal D } \abs{ \sh{\mathcal E(I,k) } }\le{} k_0\abs{\sh{\mathcal U } },  \qquad k\le{}k_0. 
 \end{equation*}
And for $k>k_0 $, we have $\mathcal E(I,k)=\emptyset $ for all $I $.   By (\ref{e.2nd}) this completes the proof.

%%%%%%%%%%%%%%%%%%%%%%%%%%%%%%%%%%%%%%%%%%%%%%%%%%%%%%%%%%%%%%
\subsection{Uniform Embeddedness in Two Parameters}\parskip=11pt 
\label{s.uniform}
  
We define the  notion of embeddedness by simultaneously expanding all sides of the rectangle. 
Let $U$ be a subset of  $\mathbb R^2$ of finite measure.  We inductively define  a sequence of enlarged sets associated to $U$ by 
	  \begin{gather} \label{e.enl}
	 \Enl 2,\mathcal U, {}{}\eqdef{}\{M \ind {\sh{\mathcal U}} > \tfrac1{16}\}, \\ \Enl {j+1},\mathcal U, {}{}\eqdef{} \Enl 2, {\text{\rm Enl}_{ j}( U)  } ,\quad  j>2.
	  \end{gather}
	 Given a dyadic rectangle $R\in\mathcal U$, we give measures of how deeply embedded this rectangle is inside of $U$ by 
	  \begin{equation} \label{e.emb}
	 \emb R. \Enl j,\mathcal U,.{}\eqdef{}\sup\{\mu\ge1\mid \mu R\subset\Enl j,\mathcal U, \}, \quad j\ge2.
	  \end{equation}
One can construct examples in which for many rectangles, this measure of embeddedness is essentially 
smaller than the measure used above.

  \begin{lemma}\label{l.journe-uniform-2}   For all $\epsilon>0$, for all  collections of rectangles $\mathcal U$,  whose shadow 
 has finite measure in the plane, we have the inequality 
  \begin{equation*}
 \sum_{R\in\mathcal U'}{\emb  R. \Enl 2,\mathcal U,.}^{-\epsilon}\abs R {}\lesssim{}\abs{\sh{\mathcal U'}}.
  \end{equation*}
 The implied constant depends only on $\epsilon$, and holds uniformly over all collections $\mathcal U'\subset\mathcal U$.  
  \end{lemma}
 
 This form of the Journ\'e Lemma was first proved in Ferguson and Lacey \cite{MR1961195}.

\medskip
\subsubsection{The First Proof.} 
We rely very much on the version of Journ\'e's Lemma that we have already established.  Indeed, we will 
need a variant of this Lemma, one in which the standard dyadic grid is replaced by a shifted dyadic grid, as 
defined in Section~\ref{s.notation}. See in particular (\ref{e.shifted-grids}).  

Apply Lemma~\ref{l.journe-classic}, to a collection of rectangles $\mathcal U $.  For an integer $k $, we 
consider a collection of  rectangles $R\in\mathcal U $ such that $\emb R.\mathcal U.\le{} 2^k $, where the embeddedness 
quantity is defined as in (\ref{e.1st-emb}).  Call this collection $\mathcal U' $. 

We now define a new collection $\mathcal V $ of rectangles.  These rectangles will be  a product of $\mathcal D_{\mathsf 1} $
and a dyadic interval. Recall that for the collection of intervals $\mathcal D_{\mathsf 1} $, for any interval
$K\subset \mathbb R $, we can find $I,I'\in\mathcal D_{\mathsf 1} $ so that 
 \begin{equation*}
\tfrac 14 I\subset K \subset I'\subset 2K . 
 \end{equation*}

Now, for each $I\times J\in \mathcal U' $, take $\widetilde I\in\mathcal D_{\mathsf 1 } $ to be the maximal element such 
that (i) $I\subset\widetilde I $ and (ii) 
 \begin{equation*}
  \widetilde I\times \tfrac{\abs{\widetilde I }}{\abs I} J \subset  \Enl 2,\mathcal U,.
 \end{equation*}
 Set $\mathcal V{}\eqdef{}\{\widetilde I\times J \mid I\times J\in\mathcal U' \} $.  
Certainly, we have by the Journ\'e Lemma,
 \begin{align*}
\abs{ \sh{\mathcal V } }& {}\lesssim{}\abs{\sh {\mathcal U } }, 
\\
\sum _{R\in\mathcal U'}\abs R &{}\lesssim{} 2 ^{\epsilon k }\abs{\sh {\mathcal U } }.
 \end{align*}
In the second line, $\epsilon >0 $ is an arbitrary positive constant, and the implied constant depends upon 
$\epsilon $. 

Clearly, we want to apply the Journ\'e Lemma to the collection $\mathcal V $ in the second coordinate. This is not 
quite straight forward to do, as the collection of rectangles $\mathcal V $ may not consist exclusively of pairwise 
incomparable rectangles.  Yet, if we have two rectangles 
$\widetilde I\times J\subset \widetilde I'\times J' $, with both rectangles in the collection $\mathcal V $,  
 and in addition we have 
 \begin{equation*}
8\abs {\widetilde I}\le{}\abs {\widetilde I' }, \qquad 
 2^{k+2}\abs J\le{}\abs {J'} . 
 \end{equation*}
then, it would be the case that $\emb I\times J.\mathcal U.>2^k $, which is a contradiction.  Therefore, we see that 
$\mathcal V $ is a union of at most $O(k) $ subcollections  $\mathcal V' $,  each  of 
which consists only of pairwise incomparable rectangles.  Thus, we deduce from Lemma~\ref{l.journe-classic} that 
 \begin{equation*}
\sum _{R\in\mathcal V' }\abs R {}\lesssim{}2 ^{\epsilon k} \abs{ \sh {\mathcal V } } {}\lesssim{}  
	2 ^{\epsilon k} \abs{ \sh {\mathcal U } }.
 \end{equation*}
Therefore, the proof is complete.

\subsubsection{The Second Proof.} 

We employ the standard reduction (\ref{e.standard}), and use the ``essentially disjoint'' argument to prove the Lemma. 

The main construction of the proof is this inductive procedure.  We construct a decomposition of 
 $\mathcal U$ into ``good'' $\mathcal G(\mathcal U)$  and ``bad'' $\mathcal B_j(\mathcal U)$ parts, with $j=1,2$.  Initialize 
  \begin{equation*}
 \operatorname{Stock}{}\eqdef{}\mathcal U,\quad \mathcal G=\emptyset,\quad \mathcal B_j=\emptyset, \ j=1,2. 
  \end{equation*}
  If $\operatorname{Stock}=\emptyset$ we return $\mathcal G(\mathcal U)=\mathcal G$,  $\mathcal B_j(\mathcal U){}\eqdef{}\mathcal B_j$, for $j=1,2$.
  
  While $\operatorname{Stock}$ is non--empty, select any $R\in\operatorname{Stock}$, and update 
   \begin{equation*}
  \operatorname{Stock}=\operatorname{Stock}-\{R\},\quad \mathcal G=\mathcal G\cup\{R\}.
   \end{equation*}
  Continuing, for $j=1,2$, while there is an $R'\in\operatorname{Stock}$ so that there are $R_1,R_2,\ldots,R_N\in\mathcal G$ such that 
  the $R_n$ are longer than $R'$ in the $j$th coordinate, and 
   \begin{equation} \label{e.nearlycover}
  \Abs{ R'\cap \bigcup_{n=1}^N R_n}>\tfrac89\abs{R'},
   \end{equation}
  update 
   \begin{equation*}
  \operatorname{Stock}=\operatorname{Stock}-\{R'\},\quad \mathcal B_j=\mathcal B_j\cup\{R'\}.
   \end{equation*}
  \smallskip

  \medskip

  By construction, the rectangles in $\mathcal G(\mathcal U)$ are essentially disjoint. 
  It suffices therefore to  argue that for $j=1,2$, we have 
   \begin{equation} \label{e.1(1)} 
  \mathcal B_j(\mathcal B_j (\mathcal U))=\emptyset.
   \end{equation}
  And it follows that  inductively applying the decomposition into good and bad parts  to each of $\mathcal B_j(\mathcal U)$ will terminate after three rounds.

  Suppose by way of contradiction, that there is an $R\in\mathcal B_1(\mathcal B_1(\mathcal U))$.  Thus, there 
  are $R_1,R_2,\ldots,R_N\in\mathcal B_1(\mathcal U)$ for which 
  each $R_n$ is longer in the first coordinate and (\ref{e.nearlycover}) holds.  Then, 
  suppose that $R_1$ has first coordinate $R_{1(1)}$ that  among all the $R_n$ is shortest 
  in the first coordinate.  Since each $R_n$ is in $\mathcal B_1(\mathcal U)$ each of these rectangles 
  are themselves nearly covered by rectangles in $\mathcal U$ that are longer in the first 
  coordinate.  By the standard reduction, these rectangles are themselves much longer than 
  $R_{1(1)}$.  Hence, we take $I$ to be the dyadic interval of length 
  $10\mu\abs{R_{1(1)}}\le\abs{I}<20\mu\abs{R_{1(1)}}$ that contains $R_{1(1)}$.  Let $J$ be 
  the second coordinate of $R$.  Then, it is necessarily the case that 
   \begin{equation*}
  \abs{I\times J \cap \sh {\mathcal U'}}\ge{}(\tfrac89)^2\abs{I\times J}.
   \end{equation*}
  But then, $\frac98(I\times J)\subset\Enl 2,\mathcal U,$.

   $I$ is much larger than $R_1$ in the first coordinate, as we have separated scales. 
  For the same reason, $J$ is much longer than than $R_1$ in the second coordinate.  Hence, 
  we see that $3\mu R_1\subset  \frac98(I\times J)$. But this contradicts the assumption 
  that $\emb R.\Enl 2,\mathcal U,.\le2\mu$, and so completes the proof of the Lemma.

%%%%%%%%%%%%%%%%%%%%%%%%%%%%%%%%%%%%%%%%%%%%

 %%%%%%%%%%%%%%%%%%%%%%%%%%%%%%%%%%%%%%%%%%%%%%%%%%%%%%%%%%%%%%
 \subsection{Uniform Embeddedness with Small Enlargement in Two Parameters}\parskip=11pt
 \label{s.uniform-small}
 
 In this section, our emphasis shifts to the enlarged sets.  Specifically, we permit  the enlarged set  $\Enl ,\mathcal U,$
	 to be only slightly bigger than $\sh{\mathcal U}$ itself, no more $\abs {\Enl ,\mathcal U,}\le{}(1+\delta)\abs{\sh{\mathcal U}}$, where $\delta>0$ is arbitrarily small. 
	 
	 We shall see that as $\delta$ decreases, the method by which we have to select it changes considerably.  
	 So let us emphasize that $U\subset V$, and that we shall define 
	  \begin{equation*}
	 \emb R.V.{}\eqdef{}\sup\{\mu\,\ge1\mid \mu R\subset V\},\qquad R\in\mathcal U.
	  \end{equation*}

The fact we wish to explain is the Lemma from the Appendix of \cite{MR1961195}.\footnote{It seems likely  one could also use the embeddedness in Section~\ref{s.redux}, but we do not pursue that here.}

 \begin{proposition}\label{p.smallV}  
  For    each $0<\delta,\epsilon <1$  there is a constant $K_{\delta,\epsilon}$, so that for all  
  for all collections of rectangles $\mathcal U$ whose shadow has finite measure in the plane, 
   there is a set $V\supset\sh{\mathcal U}$ for which $\abs V<(1+\delta)\abs{\sh{\mathcal U}}$,  so that for any collection  $\mathcal U'\subset\mathcal U$ 
  we have the inequality 
   \begin{equation} \Label e.JsmallV  \sum_{R\in\mathcal U'} \emb R.V.^{-\epsilon}\abs R\lesssim{}\abs{\sh{\mathcal U'}}. 
   \end{equation}
  The implied constant depends only  on $\epsilon,\delta>0$.
   \end{proposition} 

 We define $V$. Recall the definition and 
properties of shifted dyadic grids, (\ref{e.shifted-grids}).
 For a collection of intervals $\mathcal I$ and $j=1,2$, set $\operatorname M^\mathcal I_j$ to be
 the maximal function associated to $\mathcal I$, computed in the coordinate $j$.  Initially, we use only
 the dyadic grids,  setting  $\delta=(1+2^{\mathsf d})^{-1}$ and 
  \begin{equation*}
\Enl 0,\mathcal U,{}\eqdef{}\bigcup_{i\not=j}\{\operatorname M^\mathcal D_i \ind {} \{\operatorname M_j\ind {\sh{\mathcal U}}>1-\delta\}>1-\delta\}.
  \end{equation*}
 It is clear that $\abs {\Enl 0,\mathcal U,}<(1+K\delta)\abs{\sh{\mathcal U}}$.  Invoking the collections $\mathcal D_{\mathsf d}$, set 
  \begin{equation} \label{e.V}
 \Enl ,\mathcal U,{}\eqdef{}\bigcup_{i\not=j}\{\operatorname M_i^{\mathcal D_{\mathsf d}}\ind {} \{
 \operatorname M^{\mathcal D_{\mathsf d}}_j\ind {V_0}>1-\delta\}>1-\delta\}.
  \end{equation}
 Then $\abs {\Enl ,\mathcal U,}<(1+K\delta\log\delta^{-1})\abs{\sh{\mathcal U}}$, and we will work with this choice of $V$.  
 This is the set $V$ of the Lemma.

 The additional important property that $\Enl ,\mathcal U,$ has can be formulated this way.  For all dyadic
 rectangles $R=R_1\times R_2\subset \Enl 0,\mathcal U,$, the four rectangles 
  \begin{equation} \Label e.Dd  
 (R_1\pm\delta\abs{R_1})\times(R_2\pm\delta\abs{R_2})\subset\Enl ,\mathcal U,
  \end{equation}
 This follows immediately from the construction of  the shifted dyadic grids.  The first
 stage of the proof is complete.

 \medskip

The remainder of the argument is as in Section~\ref{s.uniform}.
We impose the standard reduction, with the additional stipulation that the scales in $\mathcal U$ be separated by $10^{6}\mu\delta^{-1}$.  
And we use the essentially disjoint proof strategy.
There is  a ``bad'' class of rectangles  $\mathcal B=\mathcal B(\mathcal U)$
to consider, defined as follows.  For $j=1,2$, let $\mathcal B_j(\mathcal U)$ be 
 those rectangles $R$ for which there are rectangles 
  \begin{equation*} 
 R^1,R^2,\ldots,R^K\in\mathcal U-\{R\}, 
  \end{equation*} 
 so that for each $1\le{}k\le{}K$, $\abs{R^k_j}>\abs{R_j}$, and 
  \begin{equation*} 
 \ABS{ R\cap \bigcup_{k=1}^K R^k}>(1-\tfrac\delta{10})\abs R. 
  \end{equation*} 
 Thus $R\in\mathcal B_j$ if it is nearly completely covered by dyadic rectangles 
 in  the $j$th direction of  the plane.  
Set $\mathcal B(\mathcal U)=\mathcal B_1(\mathcal U)\cup\mathcal B_2(\mathcal U)$.      It follows that 
if $R\not\in\mathcal B(\mathcal U)$, it is not covered in both the vertical and horizontal 
directions, hence 
 \begin{equation*} \ABS{ R\cap 
\bigcap_{R'\in\mathcal U-\{R\}}(R')^c}\ge\tfrac{\delta^2}{100}\abs R. 
 \end{equation*} 
And so  \begin{equation*} 
\sum_{R\in\mathcal U-\mathcal B(\mathcal U)}\abs R\le100\delta^{-2}\abs{\sh{\mathcal U}}. 
 \end{equation*} 

\medskip 
 
Thus, it remains to consider   the set of rectangles $\mathcal B_1(\mathcal U)$ 
and $\mathcal B_2(\mathcal U)$.      Observe that for any collection $\mathcal U'$, 
$\mathcal B_j(\mathcal U')\subset\mathcal U'$ as follows 
immediately from the definition.  Hence 
$\mathcal B_1(\mathcal B_2(\mathcal B_1(\mathcal U)))\subset{} \mathcal B_1(\mathcal B_1(\mathcal U))$.  And we argue that this last set is empty.  
As our definition of $V{}\eqdef{}\Enl ,\mathcal U,$ and $\emb R.\mathcal U.$ is symmetric with respect to the 
coordinate axes, this is enough to finish the proof. 
 
 We argue that $\mathcal B_1(\mathcal B_1(\mathcal U))$ is empty by contradiction. 
Assume that $R $ is in this collection.  Consider 
  those rectangles $R'$ in $\mathcal B_1(\mathcal U)
  $ for which $(i)$ $\abs{R'_1}>\abs{R_1}$  and 
$(ii)$ 
$R'\cap R\not=\emptyset$. 
Then  \begin{equation*} \ABS{R\cap\bigcup_{R'\in\mathcal B_1(\mathcal U)}R'}\ge(1-\tfrac\delta{10})\abs R. 
 \end{equation*} 
Fix a  one of these rectangles $R'$ with 
 $\abs{R'_1}$ being minimal.  We then claim that 
$8\mu R'\subset  \Enl ,\mathcal U ,$, which 
contradicts the assumption that $\emb R'.\mathcal U.$ is no more than $2\mu$. 
 
Indeed,  all the rectangles in $\mathcal B_1(\mathcal U)$ are themselves covered by dyadic rectangles in 
the first  coordinate axis. We see that the the set  $\{M_2^\mathcal D\ind {\sh{\mathcal U}}>1-\delta\} $ 
contains the dyadic rectangle $R''_1\times  R_2$, in which $R_2$ 
is the second coordinate interval for the rectangle $R$ and 
$R''_1$ is the dyadic interval that contains $R'_1$ and has 
measure $8\mu\delta^{-1}\abs{R_1'}\le{} 
\abs{R''_1}<16\mu\delta^{-1}\abs{R_1'}$. 
 
That is $R''_1\times  R_2$ is contained in $\Enl 0,\mathcal U,$. And the  dimensions of this rectangle are very
much bigger than those of $R$. Applying (\ref{e.Dd}), the rectangles $ (R''_1\pm\abs{R''_1})\times
R_2\pm\delta\abs{R_2}$  are  contained in  $\Enl , \mathcal U,$. And 
since $8\mu R'$ is contained in one of these  last four  rectangles, we have 
contradicted the assumption that $\emb R'.\mathcal U.<2\mu$.

 %%%%%%%%%%%%%%%%%%%%%%%%%%%%%%%%%%%%%%%%%%%%%%%%%%%%%%%%%%%%%%
 
 \subsection{Uniform Embeddedness Redux} \parskip=11pt
 \label{s.redux}
 
 The previous notion of embeddedness expanded all directions in an equal amount. We propose here  
	 an alternate method, in which non diagonal dilations are used.\footnote{Unlike the other formulations of 
	 Journ\'e's Lemma in this paper, this 
	 one has not as of yet found application in the literature.} 
	 We continue with  the definitions of (\ref{e.enl}). For 
	 a vector of positive numbers $(\mu_1,\mu_2)$, set
        \begin{equation*}
       \emb R. \mathcal U.{}\eqdef{}{} 
       \sup\{  \mu_1\mu_2 \mid \dil (\mu_1,\mu_2); R;\subset  \Enl 2,\mathcal U, ,\ 
       	\mu_1,\, \mu_2\ge1\}.
	 \end{equation*}
	This definition of embeddedness can be essentially smaller than the form studied in Section~\ref{s.uniform}. 

 \begin{lemma}\label{l.journe-inf}  In the case $d=2$, for any $\epsilon>0$, any 
collection of rectangles $\mathcal U$ in the plane, whose shadow has finite measure, 
and all $\mathcal U'\subset\mathcal U$ of rectangles which are maximal, we have 
 \begin{equation*}
\sum_{R\in\mathcal U'} \emb R.\mathcal U'.^{-\epsilon}\abs R\lesssim{}\abs{\sh{\mathcal U'}}.
 \end{equation*}
The implied constant depends only on $\epsilon>0$. 
 \end{lemma}

To prove the Journ\'e Lemma, we assume that $\mathcal U$ satisfies  the standard reduction.  We should 
further refine this reduction.  Fix  
$(\mu_1,\mu_2)$ with $\mu\le\prod_{j=1}^2 \mu_j\le2\mu$ and each $\mu_j\ge1$.  
Assume that for each $R\in\mathcal U$ we have 
 \begin{equation*}
\dil  (\mu_1/2,\mu_2/2);R;\subset  \Enl 2,\mathcal U, \qquad  \dil 2(\mu_1,\mu_2);R;\not\subset  \Enl 2, U,.
  \end{equation*}
 It suffices to consider $\lesssim(\log \mu)^3$ such collections. The argument of Section~\ref{s.uniform} proceeds with only modest changes.

%%%%%%%%%%%%%%%%%%%%%%%%%%%%%%%%%%%%%%%%%%%%%

\section{High Parameter Case with Unions of Rectangles}\parskip=11pt
\label{s.few}

We introduce a  variant of Journ\'e Lemma, in parameters three and higher, which 
can be found in J.~Pipher's paper \cite{MR88a:42019}. 
We measure embeddedness in only  one coordinate, but then must form the sum over sets more general 
than rectangles.  

The necessity of this can be seen by considering a set in $\mathbb R^3$ of the form $U=[0,1]\times{}U_2$, for a set $U_2\subset\mathbb R^2$. 
Each rectangle $R=[0,1]\times{}R_2\subset U$  has a measure of embeddedness in the first coordinate of $1$.  But the rectangles 
are certainly not disjoint in general in the second and third coordinates.

 \subsection{With Large Enlargement}
 Let $U$ be a subset of $\mathbb R^d$ with finite measure. Our Lemma makes sense in two parameters, 
but is primarily of  interest in parameters $d\ge3$.  Let $\mathcal U$ be a set of maximal 
dyadic rectangles contained in $U$.  And define the enlarged set, and embeddedness by 
 \begin{gather} \label{e.d-11} 
\Enl ,\mathcal U,{}\eqdef{}\{ \operatorname M \ind {\sh{\mathcal U}}>\tfrac12\}.
\\  \label{e.d-12}
\emb R.\mathcal U.{}\eqdef{}\sup\{\mu\ge1\mid \dil (\mu,1,\ldots,1); R;\subset \Enl ,\mathcal U,\}.
 \end{gather}
In the top line, the first maximal function $\operatorname M_1$ is applied in the first coordinate only, 
and the second maximal function $\operatorname M$ is the usual strong maximal function.\footnote{
The second maximal function could be restricted to the strong maximal function in all coordinates {\em except} the first coordinate.} 
The sets which we sum over are  specified by the choice of subcollection $\mathcal U'\subset\mathcal U$, a choice of $j\in\mathbb N$ and dyadic interval $I\in\mathcal D$. 
 \begin{equation*}
F(I,j,\mathcal U'){}\eqdef{}\bigcup\{ I\times R'\mid I\times R'\in \mathcal U',\ 2^{j-1}\le\emb  I\times R'.\mathcal U.<2^j\}.
 \end{equation*}
Our Lemma is then\footnote{We have stated the lemma in the formulation for the first 
coordinate to ease the burden of notation. In application,  the role 
of the first coordinate is  imposed on an arbitrary choice of coordinate.}

 \begin{lemma}\label{l.few} For all $\epsilon>0$, we have the estimate 
 \begin{equation*}
\sum_{j=1}^\infty \sum_{I\in\mathcal D} 2^{-\epsilon j}\abs{ F(I,j,\mathcal U')}\lesssim{}\abs{\sh {\mathcal U'}}.
 \end{equation*}
In fact, we have the estimate below, valid for any integer $n>1$, and choice of $1<p<\infty$. 
 \begin{equation*}
\NOrm \sum_{j=1}^\infty \sum_{I\in\mathcal D} 2^{-\epsilon j}(\operatorname M\ind { F(I,j,\mathcal U')})^n .p.\lesssim{}\abs{\sh {\mathcal U'}}^{1/p}.
 \end{equation*}
These estimates hold for all collections of rectangles $\mathcal U$ whose shadow has finite measure, and all collections $\mathcal U'\subset\mathcal U$.
 \end{lemma}

 This estimate has the most power when one has a collection of 
rectangles $\mathcal U$ for which the rectangles are embedded in the enlarged set by a small 
amount in the first coordinate, but embedded in the other coordinates by a very large amount. 

\bigskip 

We begin the proof, which draws upon the arguments of J.~Pipher \cite{MR88a:42019}. Fix an integer $j$. We assume that $2^{j-1}\le{}\emb R.\mathcal U.<2^j$ for all $R\in\mathcal U$, and separate scales accordingly.  Following a modification of the ``essentially disjoint'' proof strategy,  we then identify a subset $H(I)\subset F(I,j,\mathcal U')$ for which  
 \begin{equation*}
\abs{ H(I)}\gtrsim\abs{F(I,j,\mathcal U')},\qquad \qquad \operatorname M\ind {H(I)}\gtrsim\ind {F(I,j,\mathcal U')}, 
 \end{equation*}
and these sets are disjoint as $I\in\mathcal D$ varies. The first claim of the Lemma is then clear, and the second claim follows from the Fefferman--Stein maximal function estimate. 

Define 
 \begin{equation*}
G(I){}\eqdef{}\bigcup_ {I\subset_{\not=}I'}F(I',j,\mathcal U').
 \end{equation*}
Suppose that it is the case that for some $R\in\mathcal U$ with $\zR1=I$, that we have 
 \begin{equation*}
\abs{R\cap G(I)}\ge\tfrac34\abs{R}. 
 \end{equation*}
Then, by separation of scales, we see that $\emb R.\mathcal U.>2^j$, a contradiction. 
Therefore, we take the set $H(I){}\eqdef{}F(I,j,\mathcal U')\cap G(I)^c$.  These sets are clearly disjoint in $I$.  

By construction, we must have that $\abs{ R\cap H(I)}\ge\tfrac14\abs{R}$, for all $R\in\mathcal U$ with $\zR1=I$. 
Hence, applying the strong maximal function, we see that 
 \begin{equation*}
\abs{F(I,j,\mathcal U')}\le{}\abs{\{\operatorname M\ind {H(I)}\ge\tfrac14\}}\lesssim\abs{H(I)}.
 \end{equation*}
This completes the proof.

%%%%%%%%%%%%%%%%%%%%%%%
\subsection{ With Small Enlargement} \parskip=11pt

\label{s.small-few}

There is a version of the previous Lemma that employs an enlargement that is only 
slightly larger than the set $U$, as in Section~\ref{s.uniform-small}.
    
Given a   a collection of rectangles $\mathcal U$ whose shadow has finite measure,
suppose that $V\supset U$, and define 
 \begin{equation*}
\emb R.V.{}\eqdef{}\sup\{\mu\ge1\mid \dil (\mu,1,\ldots,1);R;\subset V\}.
 \end{equation*}
As before, define 
 \begin{equation*}
F(I,j,\mathcal U'){}\eqdef{}\bigcup\{ I\times R'\mid I\times R'\in \mathcal U',\ 2^{j-1}\le\emb  I\times R'.V.<2^j\}.
 \end{equation*}

 \begin{lemma}\label{l.few-small} For all $\delta,\epsilon>0$, all  $\mathcal U$ as above, we can select $V\supset \sh{\mathcal U}$ 
with $\abs V\le{}(1+\delta)\abs{\sh{\mathcal U}}$, for which   we have the estimate 
 \begin{equation*}
\sum_{j=1}^\infty \sum_{I\in\mathcal D} 2^{-\epsilon j}\abs{ F(I,j,\mathcal U')}\lesssim{}\abs{\sh {\mathcal U'}}.
 \end{equation*}
This holds for all sets $\mathcal U'\subset\mathcal U$. 
In fact, we have the estimate below, valid for any integer $n>1$, and choice of $1<p<\infty$. 
 \begin{equation*}
\NOrm \sum_{j=1}^\infty \sum_{I\in\mathcal D} 2^{-\epsilon j}(\operatorname M\ind { F(I,j,\mathcal U')})^n .p.\lesssim{}\abs{\sh {\mathcal U'}}^{1/p}.
 \end{equation*}
The implied constants in these estimates depend only on dimensions and the choices of $\epsilon,\delta$.  
 \end{lemma}

We begin the proof.  Recall the properties of shifted dyadic grids (\ref{e.shifted-grids}).  We take 
$\delta=(1+2^{\mathsf d})^{-1}$, for integer $\mathsf d$. We use the maximal function 
$\operatorname M^{\mathcal D_{\mathsf d}}$, which satisfies (\ref{e.small-weak}).  Define 
 \begin{equation*}
V{}\eqdef{}\{ \operatorname M^{\mathcal D_{\mathsf d}}_1\ind {\text{Enl}_{1} (\mathcal U)}>1-\delta\},
 \end{equation*}
 Then, it is the case that $\abs V\le{}(1+K \delta)\abs{\sh{\mathcal U}}$.

The remainder of the proof is much as in the previous section.
 We assume that $\mathcal U'$ is such that $2^j\le{}\emb R.V.<2^{j+1}$, for all $R\in\mathcal U'$,  and separate 
scales by $40\cdot 2^j $.

Define 
 \begin{equation*}
G(I){}\eqdef{}\sup_{I\subset_{\not=}I'}F(I',j,\mathcal U').
 \end{equation*}
Suppose that it is the case that for some $R\in\mathcal U'$ with $\zR1=I$, that we have 
 \begin{equation*}
\abs{R\cap G(I)}\ge(1-\tfrac\delta2)\abs{R}. 
 \end{equation*}
Then, by separation of scales, we see that $\emb R.\mathcal U.>2^j$, a contradiction. 
Therefore, we take the set $H(I){}\eqdef{}F(I,j,\mathcal U')\cap G(I)^c$.  These sets are clearly disjoint in $I$.  

By construction, we must have that $\abs{ R\cap H(I)}\ge\frac\delta2\abs{R}$, for all $R\in\mathcal U$ with $\zR1=I$.  Hence, applying the strong maximal function, we see that 
 \begin{equation*}
\abs{F(I,j,\mathcal U')}\le{}\abs{\{\operatorname M\ind {H(I)}\ge\tfrac\delta 4\}}\lesssim\delta^{-1}\abs{H(I)}.
 \end{equation*}
The construction of these sets proves the Lemma.  

%%%%%%%%%%%%%%%%%%%%%%%

\subsection{With Uniform Embeddedness} 

We list a  version of Lemma~\ref{l.few} which has some advantages as we   use a uniform notion of embeddedness. 

 \begin{lemma}\label{l.few-uniform}  For all $\epsilon>0$, and all collections of rectangles $\mathcal U$ whose shadow has 
finite measure, there is a set $V\supset \sh {\mathcal U} $ such that $\abs V{}\lesssim{} \abs{\sh {\mathcal U}} $, 
and there is a map $\operatorname {emb} \mid \mathcal U \mapsto [1,\infty)$, and a map
$\imath \mid \mathcal U\mapsto \{1,2,\ldots,d\} $  such that 
 \begin{equation*}
\operatorname {emb}(R)\cdot R\subset V,\qquad R\in\mathcal U 
 \end{equation*}
and for all collections $\mathcal U'\subset\mathcal U$, 
 \begin{gather*}
\sum_{j=1}^d\sum_{v=0}^\infty 
\sum_{I\in\mathcal D}  2^{-(d+\epsilon) v}\abs{ F(I,j,v,\mathcal U')}\lesssim{} \abs{\sh{\mathcal U'}} ,
\\
\text{where }  \qquad F(I,j,v,\mathcal U'){}\eqdef{}\bigcup\{R\in\mathcal U'\mid  2^v<\operatorname {emb}(R)\le2^{v+1},\  \zR j=I,\ \imath (R)=j\}.
 \end{gather*} 
 \end{lemma} 

Notice that we have to have a substantially worse power on the embeddedness term, namely the power of embeddedness is 
strictly smaller than $-d $.

The method of proof requires that we  apply Lemma~\ref{l.few-small}, although
we find it necessary to apply it both inductively and to a wide range of possible collections of
rectangles.  In fact, it is useful to us that Lemma~\ref{l.few-small} applies not just to collections of
dyadic rectangles  $\mathcal U$ such that the shadow of $\mathcal U$ is of finite measure.   
It also applies to collections of rectangles $\mathcal U\subset\otimes _{j=1}^d\mathcal D _{\mathsf1}  $, where we are referring to the 
$\mathsf d $ fold product of shifted dyadic intervals.    And it moreover applies to all subcollections of $\mathcal U $.

We apply Lemma~\ref{l.few-small}  to $\mathcal U^0{}\eqdef{}\mathcal U$.  Thus, we get a set $V^{1}\supset\sh {\mathcal U^0}$, with $\abs{ V^1} {}\lesssim{}
\abs {\sh {\mathcal U^0}}$, so that for
 \begin{equation*}
\operatorname{emb}^1 (R,V^1){}\eqdef{}\sup\{\mu\ge1\mid \mu \zR 1\times\zR 2\times \cdots\times \zR n\subset V^1\}.
 \end{equation*}
we have the conclusion of Lemma~\ref{l.few-small} holding.
We then construct $\mathcal U^1\subset\otimes_{j=1}^d \mathcal D_{\mathsf 1}$. Set
 \begin{align*}
\mathcal U^1{}\eqdef{}\{ \Gamma\times \otimes_{j=2}^n \zR j\mid& R\in\mathcal U,\
    \Gamma\in\mathcal S,\\ {}&\qquad (\zR 1\cup \tfrac14\operatorname{emb}^1 (R,V^1)\zR 1)\subset\Gamma\subset
    \operatorname{emb}^1 (R,V^1)\zR 1\}.
   \end{align*}
Notice that we are relying on the structure of the shifted dyadic grids in this definition.

The inductive stage of the construction is this.  For $2\le m\le n$, given $\mathcal U^{m-1}\subset\otimes_{j=1}^d \mathcal D_{\mathsf 1}$,   
we apply Lemma~\ref{l.few-small} to get a set   $V^m$ satisfying
 \begin{equation*}
V^m\supset \sh{\mathcal U^{m-1}},\qquad \abs {V^m} {}\lesssim{}{}
\abs{   \sh{\mathcal U^{m-1}}}.
 \end{equation*}
The embedding function for rectangles $R\in\mathcal U^{m-1}$  is
  \begin{align*}
 \operatorname{emb}^m (R,V^{m}){}\eqdef{}\sup\{\mu\ge1\mid &
 R_{1}\times\cdots\times{}R_{m-1}\times\mu R_{m}
 \\&\qquad\times\zR {m+1}\times\cdots\times\zR n\subset{}V^{m}\}.
  \end{align*}
 And the conclusion of Lemma~\ref{l.few-small} holds.
 The collection $\mathcal U^m$ is then taken to consist of all rectangles of the form
  \begin{equation*}
 \otimes_{j=1}^{m-1} \zR j\times \Gamma\times \otimes_{j=m+2}^n\zR j
  \end{equation*}
 where $R\in\mathcal U^{m-1}$ and $\Gamma\in\mathcal S$  satisfies
  \begin{equation*}
( \zR m\cup  \tfrac14\operatorname{emb}^m (R,V^{m}) \zR m)\subset\Gamma\subset \operatorname{emb}^m (R,V^{m})\zR m.
  \end{equation*}
% This is done for all $1\le{}m\le{}n$.

 To prove our Lemma, we take $V{}\eqdef{}V^n$.   It is the case that
  \begin{align*}
 \abs{V^{n}}\lesssim{}&\abs{\sh {\mathcal U^{n-1}}}
 \\{}\lesssim{}&\abs{V^{n-1}}
 \\{}\lesssim{}&\abs{\sh {\mathcal U}}.
 \end{align*}

The definition of the embedding function is not so straight forward.  It is taken to be
 \begin{equation*}
\operatorname{emb}(R)=\tfrac1{16}\inf_{1\le{}m\le{}n}\beta^m(R)
 \end{equation*}
where $\beta^m(\cdot)$ are inductively defined below. The function $\imath(R)$ is taken to be the coordinate in which the infimum for the embedding function is achieved.

Set $\beta^1(R){}\eqdef{}\operatorname{emb}^1 (R,V^{1})$.
In the inductive step,  for $2\le m\le n$, set $
\gamma _m(R){}\eqdef{}\inf_{j<m}\beta^j(R)$.  For $1<\gamma <\gamma _m(R)$, let
 \begin{equation*}
\beta^m_\gamma (R){}\eqdef{}\operatorname{emb}^m (\varphi^m_\gamma (R),V^{m})
 \end{equation*}
where $\varphi^m_\gamma (R)\in\mathcal U^{m-1}$ is the rectangle with $\varphi^m_\gamma (R)_{j}=\zR j$ for $j\ge{}m$, and for
$1\le{}j<m$, $\varphi^m_\gamma (R)_{(j)}$ is the element of $\mathcal D _{\mathsf 1}$ of maximal length such that
 \begin{equation*}
(\zR j\cup\tfrac14\gamma \zR j)\subset \varphi^m(R)_{(j)}\subset \gamma \zR j.
 \end{equation*}
Now,  take $\overline\gamma $ to be the largest value of $1\le{}\gamma \le{}\gamma _m(R)$ for which we have the inequality
$
\beta^m_\gamma (R)\ge\gamma 
$.  Let us see that this definition of $\overline \gamma $ makes sense.
This last inequality is strict  for $\gamma =1$,  and as $\gamma $ increases, $\beta^m_\gamma (R)$ decreases, so $\overline\gamma $
is a well defined quantity.
Then define $\beta^m(R){}\eqdef{}\beta^m_{\overline \gamma }(R)$, and for our use below,
set $\varphi^m(R){}\eqdef{}\varphi^m_{\overline\gamma }(R)$.

The choices above prove our lemma, as we show now.
For each rectangle $R\in\mathcal U$, it is clear that $\operatorname{emb}(R)R\subset V$.
Take $\mathcal U'\subset \mathcal U$.
Considering  the sets $F(I,k,m,\mathcal U')$, then, by Lemma~\ref{l.few-small} applied in the $m$th coordinate,
 \begin{equation*}
\sum_{I\in\mathcal D}\abs{F(I,k,m,\mathcal U')}\le{}
2^{\epsilon k}\abs{\sh {\varphi^m(\mathcal U')}}.
 \end{equation*}
While we have a very good estimate for the shadow of $\varphi^m(\mathcal U)$, a corresponding good estimate for
an arbitrary subset $\mathcal U'$ seems very difficult to obtain.  But  it is a consequence of our construction
  that the rectangle $\varphi^m(R)$ is a rectangle which agrees with $R$ in the coordinates $j\ge{}m$ and, for
coordinates $1\le{}j<m$, is expanded by at most $32\operatorname{emb}(R)\le{}2^{k+6}$.  Hence, we have the estimate
 \begin{equation*}
\ABs{\bigcup\{ \varphi^m(R)\mid R\in \mathcal U'\}}\lesssim{}2^{d k}\abs{\sh {\mathcal U'}}.
 \end{equation*}
This follows from the weak $L^1$ bound for the maximal function in one dimension, applied in each coordinate separately.
It is in this last step that we lose the large power of the embeddedness.  Our proof is complete.

%%%%%%%%%%%%%%%%%%%%%%%%  
\subsection{With Small Enlargement}  

Continuing in this theme, there is a version of the previous Lemma in which one 
does not permit the enlarged set to be very big.  We record 

 \begin{lemma}\label{l.few-uniform-small} For all $\delta>0 $,   and  all $\epsilon>0$, there is a constant $K_ {\delta,  \epsilon } $
so that for all collections of rectangles $\mathcal U$ whose shadow has 
finite measure, there is a set $V\supset \sh {\mathcal U} $ such that $\abs V\le{}  (1+\delta)\abs{\sh {\mathcal U} } $, 
and there is a map $\operatorname {emb} \mid \mathcal U \mapsto [1,\infty)$, and a map
$\imath \mid \mathcal U\mapsto \{1,2,\ldots,d\} $  such that 
 \begin{equation*}
\operatorname {emb}(R)\cdot R\subset V,\qquad R\in\mathcal U 
 \end{equation*}
and for all collections $\mathcal U'\subset\mathcal U$, 
 \begin{gather*}
\sum_{j=1}^d\sum_{v=0}^\infty 
\sum_{I\in\mathcal D}  2^{-(d+\epsilon) v}\abs{ F(I,j,v,\mathcal U')}\le{}K_ {\delta,  \epsilon }  \abs{\sh{\mathcal U'}} ,
\\
\text{where }  \qquad F(I,j,v,\mathcal U'){}\eqdef{}\bigcup\{R\in\mathcal U'\mid  2^v<\operatorname {emb}(R)\le2^{v+1},\  \zR j=I,\ \imath (R)=j\}.
 \end{gather*} 
 \end{lemma}

This has been applied in a paper of Lacey and Terwilleger \cite{witherin}, and we refer to that paper for the detailed proof.

%%%%%%%%%%%%%%%%%%%%%%%%%%%%%%%%%%%%%%%%%%%%%%
\section{The Higher Parameter Case, with Rectangles}

The version of Journ\'e's Lemma described by J.~Pipher \cite{MR88a:42019} requires a different 
notation.  As before, we set $\Enl ,\mathcal U,{}\eqdef{}\{M\ind {\sh{\mathcal U}}>\tfrac1{2d}\}$.   For each integer $1\le{}j\le{}d$, we set 
 \begin{equation}  \label{e.eachcoordinate}
\emb j.R.{}\eqdef{}\sup\{ \mu\ge1\mid R_{(1)}\times\cdots\times  \mu R_{(j)}\times 
\cdots\times R_{(d)}\subset\Enl ,\mathcal U,\}
 \end{equation}
That is, only the $j$th coordinate of $R$ is expanded. 
While we have defined this for all coordinates $j$, we only use it for $1\le{}j<d$.  

 \begin{lemma}\label{l.journe-rectangles}  For each $d\ge3$, and $0<\epsilon<1$,   all subset $U$ of $\mathbb R^d$ of finite measure, 
and collections $\mathcal U$ of pairwise incomparable dyadic rectangles $R\in\mathcal U$, we have 
 \begin{equation*}
\sum_{R\in\mathcal U'}\abs R  \prod_{j=1}^{d-1} \emb j. R.^{-\epsilon} {}\lesssim{}\abs{\sh {\mathcal U'} }.
 \end{equation*}
The inequality holds uniformly over all subsets $\mathcal U'$ of $\mathcal U$. 
 \end{lemma}

This formulation has the advantage that the sum on the left hand side is over simpler objects, 
namely rectangles.  On the other hand, the embeddedness is now more complicated, in that it is 
a product of terms.  In particular $\abs R$ is essentially weighted by the \emph{largest} embeddedness 
term.

We give two proofs of this result.

\subsubsection{The First Proof.}
Let us define partial orders $<_j$ on dyadic rectangles by 
writing $R<_j R'$ iff $R\cap R'\not=\emptyset$ and $R_{(j)}\subset_{\not=} R_{(j)}'$.

The notion of the standard reduction is slightly different.  Let us assume that for $1\le \mu_j$, $1\le j<d$, we have a collection of 
rectangles $\mathcal U$ as above, with $\mu_j\le\emb j.R.\le2\mu_j$ for $1\le{}j<d$.  
In addition, we assume that the scales of $\mathcal U$ are separated by $10\max_j \mu_j$.  We then follow the essentially disjoint proof   
stat egy.

Suppose that there is a rectangle $R\in\mathcal U$ so that 
 \begin{equation*}
\abs{ R\cap \sh {\mathcal U-\{R\} }}\ge\tfrac78\abs R .
 \end{equation*}
Then, for some $1\le{}j\le{}d$, we can choose $R^1,\ldots,R^K\in\mathcal U-\{R\}$ with $R<_jR^k$ for all $k$, and 
 \begin{equation*}
\ABs{ R\cap \bigcup_{j=1}^KR^k }\ge\tfrac7{8d}\abs R .
 \end{equation*}
This is a contradiction to $\emb j.R.\simeq\mu_j$.   Thus, the rectangles in $\mathcal U$ are essentially disjoint, and the proof is 
complete.

\subsubsection{The Second Proof.} We give the proof of 
J.~Pipher \cite{MR88a:42019}.  It is convenient for us to restrict attention to the three parameter case, 
and comment on the higher parameter case briefly. 

We can assume that ${}\mathcal U' ${} is a collection of rectangles with ${}\emb 1.R.\simeq2^k ${} for some integer ${}k ${}.  We modify slightly the notation  of Lemma~\ref{l.few}.  
 \begin{equation*}
F(I){}\eqdef{}\bigcup\{ R\in\mathcal U' \mid R_{(1)}=I \}.
 \end{equation*}
 An essential point to observe is that if we hold the first coordinate of the rectangles ${}R ${} fixed, then 
 the two parameter arguments will apply, in particular Lemma~\ref{l.journe-classic} applies.  
 Doing so, will place Lemma~\ref{l.few} at our disposal. 
 \begin{align*}
\sum_{R\in\mathcal U'}\abs R   \emb 2. R.^{-\epsilon} &{}\lesssim{} 
	\sum_{I\in\mathcal D} \abs{ F(I)}
\\
	& {}\lesssim{}2 ^{\epsilon k} \abs{\sh {\mathcal U'}}.
 \end{align*}
This completes the proof in three parameters.  

\medskip 

In higher parameters, one can implement this proof, but one needs  certain variants of Journ\'e's Lemma 
that fall between the original formulation and Lemma~\ref{l.few}.  We omit the details.

%%%%%%%%%%%%%%%%%%%%%%%%%%%%%%%%%%%%%%%%%%%%%%%%%%%%%%%%%%%%%%

\end{document}